\documentclass[11pt,letterpaper]{article}
\usepackage[hmargin=1in, vmargin=1in]{geometry}
\usepackage{graphicx}
\usepackage{amsmath}
\usepackage{dsfont}
\usepackage{indentfirst}
\usepackage{sidecap}
\usepackage{float}
\usepackage{multicol}
\usepackage{caption}
\usepackage{subcaption}
\usepackage[hyphens]{url}
\usepackage{hyperref}

\usepackage{fancyhdr}
\title{3D Printing the Big Letters in the JDRF Logo}
\author{Edward Aboufadel \\ Grand Valley State University}
\date{Version 1.0, August 2014}

\begin{document}
\maketitle

\section{Introduction}

The purpose of this short paper is to describe a project to manufacture a 3D-print of the big letters in the logo for JDRF (formerly known as the Juvenile Diabetes Research Foundation).  The methods described in this paper involve processing the image of the logo through a Mathematica script, and the methods can be applied to most logos and other images that are sparse in the sense of not having many edges within the image.  With the Mathematica script, a STereoLithography (.stl) file is created that can be used by a 3D printer.  Finally, the object is created on a 3D printer.   We assume that the reader is familiar with the basics of 3D printing.

\section{The JDRF Logo}

JDRF (formerly known as the Juvenile Diabetes Research Foundation) is a major charity devoted to preventing, treating, and curing type-1 diabetes.  Type-1 diabetes (also called T1D) occurs when the immune system is triggered to start destroying insulin-producing beta cells in the pancreas.  The cause of this trigger is uncertain, and a cure is not in sight at this time.  People with T1D introduce artificial insulin into their bodies through injections or an insulin pump, and continuous glucose monitors are now available in order to track glucose levels in the blood at five-minute intervals.  Type-1 diabetes should not be confused with the more common type-2 diabetes that afflicts a considerable number of Americans. About 1\% of Americans have T1D.  (There is a member of the author's family afflicted with T1D.)

\begin{figure}[h]
\includegraphics[width=\textwidth]{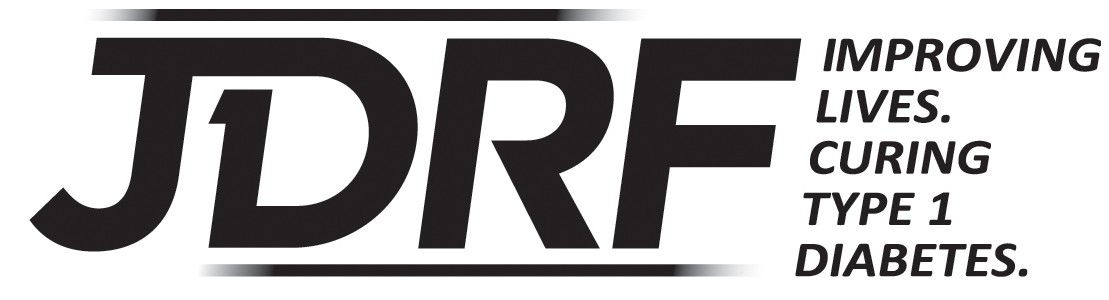}
\caption*{Figure 1: Logo of JDRF in grayscale.}
\end{figure}

A grayscale version of the logo for JDRF can be found in Figure 1. (The logo of JDRF is often seen in a specific shade of blue.)  The method described below will not work very well on the ``Improving Lives.  Curing Type 1 Diabetes'' part of the logo, at least for a 3D print the size of a hand, because of the relatively small size of the letters, and the number of letters.  However, the big ``JDRF'' is a good candidate for our method because of the relative simplicity of the four letters.  There is a also a way to address the lines which have a gradient of greyscale, located above and below the four letters.  We will describe our method in the next section.

\section{Mathematica Code}

The code in this section is modified from code created in 2013 by Melissa Sherman-Bennett and Sylvanna Krawczyk, who were my research students at our REU program.  See the URL below for more 3D printing results.

The first four lines of code loads the logo into a matrix \texttt{image} of grayscale values (on a 0 to 1 scale).  Specifically, the first line clears the memory in Mathematica, and the next lines load and process the image file from the directory \texttt{C:$\backslash$data$\backslash$3d}.  We have the \texttt{ColorConvert} command because when the Import command creates three matrices (for colors red, blue, and green), even if the jpeg file is greyscale.

\begin{verbatim}
ClearAll["Global`*"]
input = Import["C:\\data\\3d\\JDRF.jpg"];
size = Import["C:\\data\\3d\\JDRF.jpg", "ImageSize"];
image = ColorConvert[Image[input, "Real"], "Grayscale"];
\end{verbatim}

Next, we define a function that we are going to apply to every entry in \texttt{image}.  The basic idea is that the interior of the four letters will built up to a height of 5.2 millimeters from the build plate of the 3D printer.  The background of the logo will be at a height of 1.2 millimeters.  For the gradient lines, we define a gradient height, starting at 5.2 millimeters and decreasing in a reasonable way.  The way to read the following function is first there is a test that we are in a background area of the logo ($x > 0.9$ -- nearly white).  If that test fails, then we go to another test to see if we are in the area of one of the four letters ($x < 0.25$ -- very dark).  If that test also fails, we must be on one of the gradient lines, and we use the function $x -> -0.5*x + 1.3$.  Later, we will multiply the heights by 4.  Here is the function:

\begin{verbatim}
bound[x_] := If[x > 0.9, 0.3, If[x < 0.25, 1.3, -0.5*x + 1.3]];
\end{verbatim}

The next command in Mathematica is to apply this function to the image and convert to a matrix (or table).  The defined variation of $i$ and $j$ in this command takes into account that the image is stored in memory backwards from what we want to print.  There is a multiplication by 4 in this command, but that could also be implemented in the function above.

\begin{verbatim}
data = ArrayPad[Table[4*bound[ImageData[image][[i, j]]],
   {i, 1, size[[2]], 1}, {j,  size[[1]], 1, -1}], {1, 1}, 0];
\end{verbatim}

Finally, we generate an image in Mathematica, and then export that image as an STL file.  Some key options here are the \texttt{DataRange} which indicates that the print will be 80 millimeters by 28 millimeters (and as we have seen previously, 5.2 millimeters tall), and \texttt{Axes -> True, BoxRatios -> Automatic} which instructs Mathematica to create an image that is true to the dimensions defined.  The image is then exported as an STL file which can be then used on a 3D printer.

\begin{verbatim}
object = ListPlot3D[data, DataRange -> {{0, 80}, {0, 28}},
  Axes -> True, BoxRatios -> Automatic, ColorFunction -> Function[{x, y, z}, Hue[z]]]
Export["C:\\data\\3d\\JDRF.stl", object]
\end{verbatim}

\begin{figure}[h]
\includegraphics[width=\textwidth]{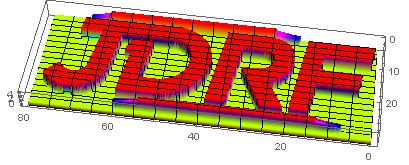}
\caption*{Figure 2: Plot of height function for our print.}
\end{figure}

The image created by the \texttt{ListPlot3D} can be found in Figure 2.  The 3D print of the letters and the gradient lines can be found in Figure 3. Another example of a logo printed using this method (and some additional tweaks) is the logo celebrating the 100th anniversary of Pi Mu Epsilon, the national honor society in mathematics.  A photo of that 3D print can be found in Figure 4.

\begin{figure}[h]
\includegraphics[width=\textwidth]{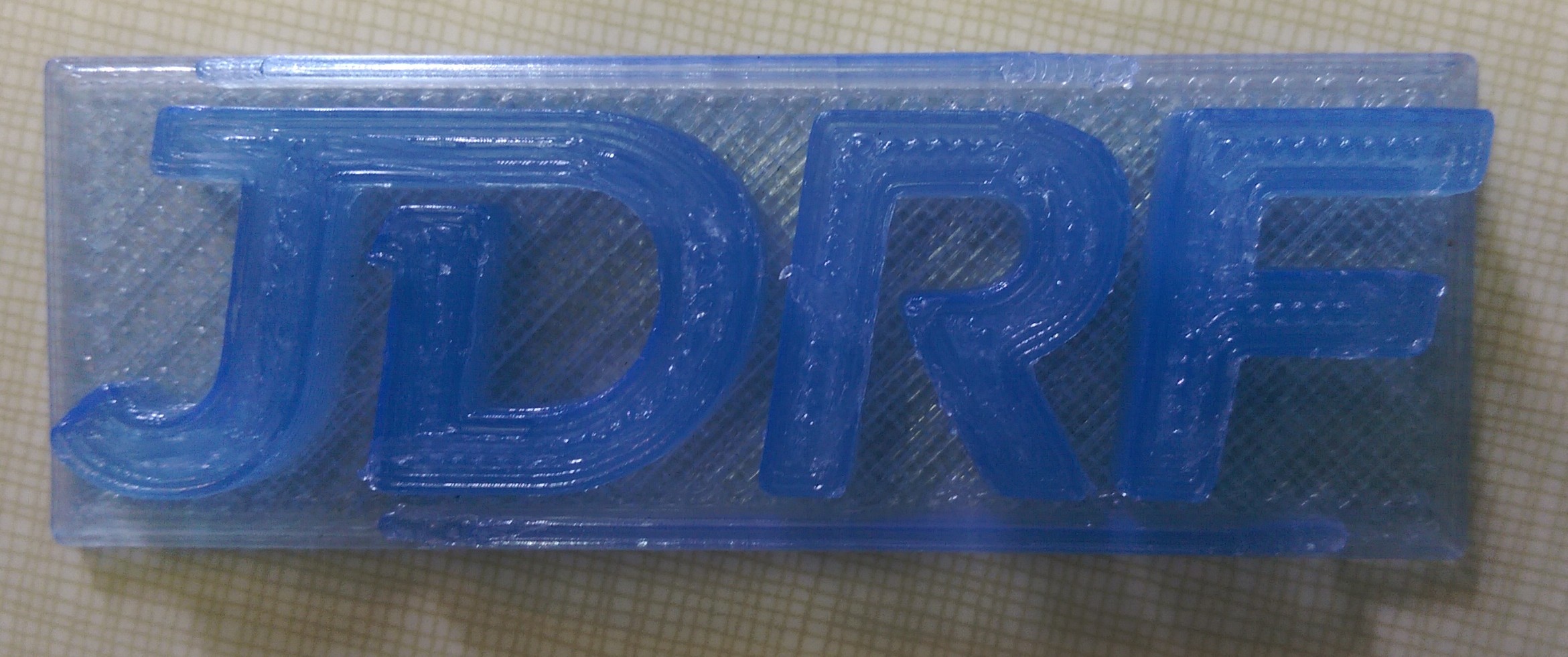}
\caption*{Figure 3: 3D Print of the Big Letters in the Logo.}
\end{figure}

\section{Appendix: URLs}
\begin{itemize}
\item JDRF home page \url{http://www.jdrf.org/}
\item Prof.~Aboufadel's 3D printing page \url{http://sites.google.com/site/aboufadelreu/Profile/3d-printing}
\end{itemize}

\begin{figure}[h]
\centering
\includegraphics[width=0.5\textwidth]{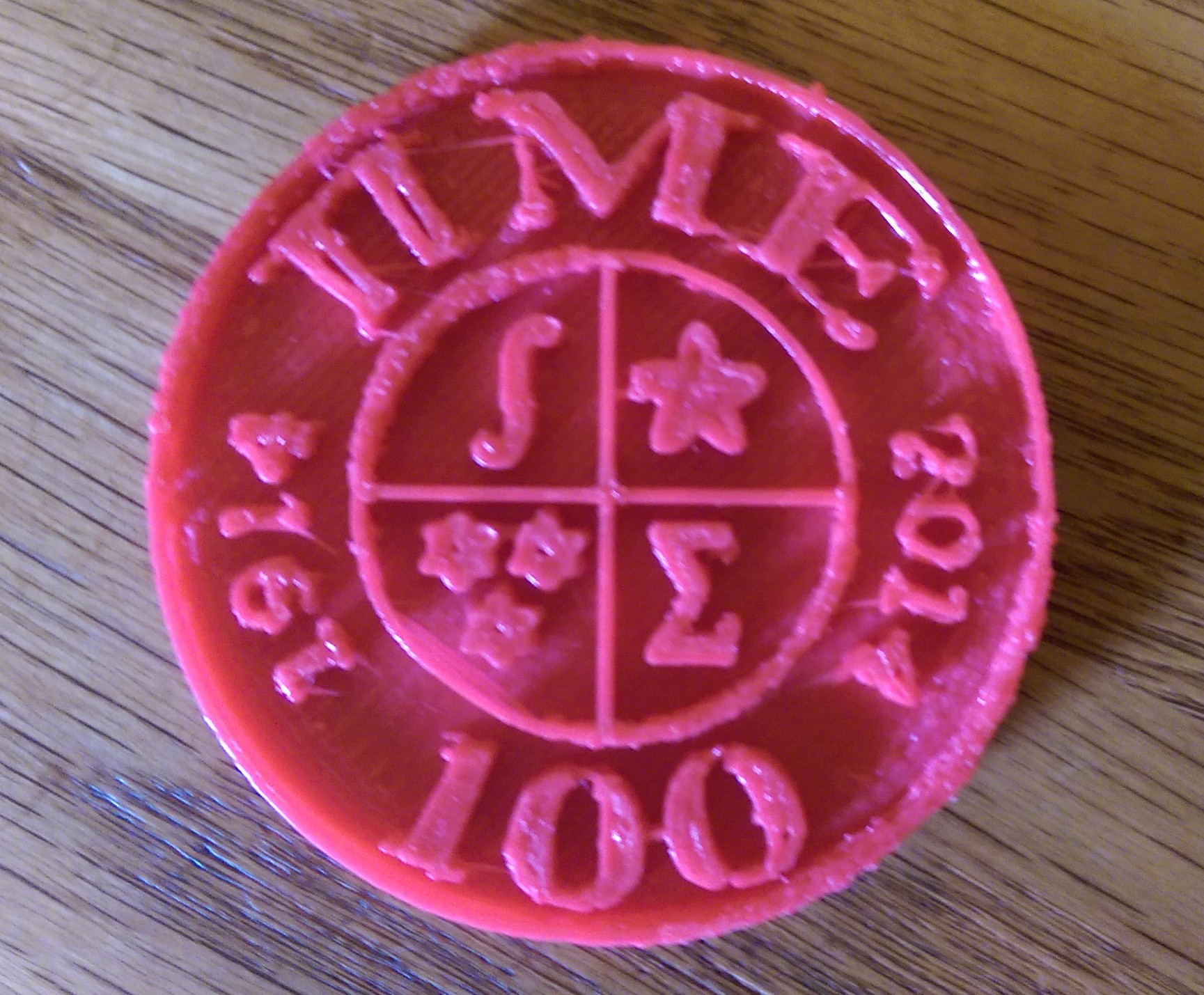}
\caption*{Figure 4: 3D Print of the logo for the PME 100th anniversary.}
\end{figure}

\end{document}